\documentclass[11pt]{article}

\usepackage{amscd,amsmath, amssymb, fancyhdr, graphicx, url,dutchcal}

\usepackage[backref=page]{hyperref}
\renewcommand*{\backref}[1]{}
\renewcommand*{\backrefalt}[4]{%
	\ifcase #1 (Not cited.)%
	\or        (Cited on page~#2.)%
	\else      (Cited on pages~#2.)%
	\fi}

\hypersetup{
	colorlinks   =true,
	citecolor    =[rgb]{0.3,0,0.3},
	linkcolor    =blue,
	urlcolor     =[rgb]{0.3,0,0.3}
}

\numberwithin{equation}{section}

\newcommand{\version}{version 2.5,\ \ June 2, 2021}

\renewcommand{\leftmark}{{\scriptsize Contraction centers in families of hyperk\"ahler manifolds}}

\makeatletter
\def\x@arrow{\DOTSB\Relbar}
\def\xlongrightarrowfill@{\arrowfill@\relbar\relbar\longrightarrow}
\newcommand{\xlongrightarrow}[2][]{%
        \ext@arrow 0099\xlongrightarrowfill@{#1}{#2}}
\makeatother

\def\eqref#1{(\ref{#1})}
\newcommand{\goth}{\mathfrak}	

\newcommand{\arrow}{{\:\longrightarrow\:}}
\newcommand{\Z}{{\Bbb Z}}
\def\C{{\Bbb C}}
\def\P{{\Bbb P}}

\newcommand{\R}{{\Bbb R}}
\newcommand{\Q}{{\Bbb Q}}

\def\1{\sqrt{-1}\:}
\newcommand{\restrict}[1]{{\left|_{{\phantom{|}\!\!}_{#1}}\right.}}
\newcommand{\cntrct}                
{\hspace{2pt}\raisebox{1pt}{\text{$\lrcorner$}}\hspace{2pt}}

\newcommand{\calo}{{\cal O}}

\renewcommand{\tilde}{\widetilde}
\renewcommand{\bar}{\overline}
\renewcommand{\phi}{\varphi}
\renewcommand{\epsilon}{\varepsilon}
\renewcommand{\geq}{\geqslant}
\renewcommand{\leq}{\leqslant}

\renewcommand{\min}{{\operatorname{\sf min}}}
\newcommand{\Teich}{\operatorname{\sf Teich}}

\newcommand{\Comp}{\operatorname{\sf Comp}}

\newcommand{\Per}{\operatorname{\sf Per}}
\newcommand{\Perspace}{\operatorname{{\Bbb P}\sf er}}

\newcommand{\Kah}{\operatorname{Kah}}

\newcommand{\Gr}{\operatorname{Gr}}

\newcommand{\Homcal}{\operatorname{{\cal H\!om}}}

\newcommand{\Pos}{\operatorname{Pos}}

\newcommand{\Aut}{\operatorname{Aut}}

\newcommand{\Mon}{\operatorname{\sf Mon}}

\newcommand{\Ext}{\operatorname{Ext}}

\newcommand{\Diff}{\operatorname{\sf Diff}}

\newcommand{\Def}{\operatorname{Def}}

\newcommand{\proof}{\noindent{\bf Proof:\ }}
\newcommand{\pstep}{{\bf Proof. Step 1:\ }}

\newcounter{Mycounter}[section]
\newcounter{lemma}[section]
\renewcommand{\thelemma}{{Lemma \thesection.\arabic{lemma}}}
\newcommand{\lemma}{%
    \setcounter{lemma}{\value{Mycounter}}
    \refstepcounter{lemma}
    \stepcounter{Mycounter}
    {\noindent \bf \thelemma:\ }}

\newcounter{claim}[section]
\newcounter{sublemma}[section]
\newcounter{corollary}[section]
\renewcommand{\thecorollary}{{Corollary \thesection.\arabic{corollary}}}
\newcommand{\corollary}{%
    \setcounter{corollary}{\value{Mycounter}}
    \refstepcounter{corollary}
    \stepcounter{Mycounter}
    {\noindent \bf \thecorollary:\ }}

\newcounter{theorem}[section]
\renewcommand{\thetheorem}{{Theorem \thesection.\arabic{theorem}}}
\newcommand{\theorem}{%
    \setcounter{theorem}{\value{Mycounter}}
    \refstepcounter{theorem}
    \stepcounter{Mycounter}
    {\noindent \bf \thetheorem:\ }}

\newcounter{conjecture}[section]
\newcounter{proposition}[section]
\renewcommand{\theproposition}
      {{Proposition \thesection.\arabic{proposition}}}
\newcommand{\proposition}{%
    \setcounter{proposition}{\value{Mycounter}}
    \refstepcounter{proposition}
    \stepcounter{Mycounter}
    {\noindent \bf \theproposition:\ }}

\newcounter{definition}[section]
\renewcommand{\thedefinition}
      {{Definition~\thesection.\arabic{definition}}}
\newcommand{\definition}{%
    \setcounter{definition}{\value{Mycounter}}
    \refstepcounter{definition}
    \stepcounter{Mycounter}
    {\noindent \bf \thedefinition:\ }}

\newcounter{example}[section]
\newcounter{remark}[section]
\renewcommand{\theremark}{{Remark \thesection.\arabic{remark}}}
\newcommand{\remark}{%
    \setcounter{remark}{\value{Mycounter}}
    \refstepcounter{remark}
    \stepcounter{Mycounter}
    {\noindent \bf \theremark:\ }}

\newcounter{problem}[section]
\newcounter{question}[section]
\makeatletter

\@addtoreset{equation}{section} \@addtoreset{footnote}{section}
\makeatother

\def\blacksquare{\hbox{\vrule width 5pt height 5pt depth 0pt}}
\def\endproof{\blacksquare}

\begin{document}

\begin{center}
{\Large\bf
Contraction centers \\[2mm] in families of hyperk\"ahler manifolds\\[3mm]
}

Ekaterina Amerik\footnote{Partially supported 
by  the  Russian Academic Excellence Project '5-100' and by French-Brasilian Research Network.}, 
Misha Verbitsky\footnote{Partially supported 
by  the  Russian Academic Excellence Project '5-100',
FAPERJ E-26/202.912/2018 and CNPq - Process 313608/2017-2.

{\bf Keywords:} hyperk\"ahler manifold, K\"ahler cone, hyperbolic geometry, cusp points

{\bf 2010 Mathematics Subject
Classification:} 53C26, 32G13}

\end{center}

{\small \hspace{0.15\linewidth}
\begin{minipage}[t]{0.7\linewidth}
{\bf Abstract} \\
We study the exceptional loci of birational
(bimeromorphic) contractions of a hyperk\"ahler manifold $M$.
Such a contraction locus is the union of all
minimal rational curves in a collection
of cohomology classes which are orthogonal to a
wall of the K\"ahler cone. Homology classes which 
can possibly be orthogonal to a wall of the K\"ahler 
cone of some deformation of $M$ are called 
{\bf MBM classes}. We prove that
all MBM classes of type (1,1) can be represented 
by rational curves, called {\bf MBM curves}. Any 
MBM curve can be contracted on an appropriate birational 
model of $M$, unless $b_2(M) \leq 5$. When $b_2(M)>5$, 
this property can be used as an alternative definition of an MBM
class and an MBM curve. Using the results of Bakker and Lehn, we prove that 
the stratified diffeomorphism type of a contraction locus 
remains stable under all deformations for which
these classes remains of type (1,1), unless the 
contracted variety has $b_2\leq 4$. Moreover, 
these diffeomorphisms preserve the MBM curves,
and induce biholomorphic maps on the contraction 
fibers, if they are normal. 
\end{minipage}
}

\tableofcontents


\section{Introduction}


\subsection{Teichm\"uller spaces in hyperk\"ahler
  geometry}
\label{_Teich_hk_Subsection_}

Let $M$ be a complex manifold. Recall that 
{\bf the Teichm\"uller space} $\Teich$ of complex structures on $M$
is the quotient $\Teich:=\Comp/\Diff_0$, where $\Comp$ is the space
of complex structures (with the topology of uniform convergence
of all derivatives) and $\Diff_0$ the connected component of the unity 
of the diffeomorphism group. The mapping class group $\Diff/\Diff_0$
on $\Teich$. 

In our case $M$ is a compact K\"ahler irreducible holomorphically symplectic
manifold \footnote{By the Calabi-Yau theorem, this is the
  same as a hyperk\"ahler manifold.} (IHS).
We consider the Teichm\"uller space $\Teich$ of all complex structures of
hyperk\"ahler type (Subsection \ref{_Teich_subsection_}). Up to the action of the subgroup $K$ of the mapping class group acting trivially on cohomology, this is the same space as the moduli space of marked
hyperk\"ahler manifolds often considered in algebraic geometry. The subgroup $K$ permutes components of $\Teich$, so that we
shall often make no difference between the relevant component of $\Teich$ and of the moduli space of marked
hyperk\"ahler manifolds\footnote{The latter has only finitely many connected components, see \cite{_Huybrechts:finiteness_}}. 

We always consider the component containing the 
parameter point of our given complex structure, i.e. parametrizing hyperk\"ahler deformations of $M$. By abuse of notation,
this space is
also denoted $\Teich$. The action of the mapping class group
$\Gamma$ (i.e. the subgroup of $\Diff/\Diff_0$ preserving our connected component) on $\Teich$ is ergodic,
and its orbits are classified using Ratner's 
orbit classification theorem (\ref{orbits-teich}).

In \cite{_Markman:universal_}, Theorem 1.1, E. Markman
has constructed the universal family on the marked moduli space. Let 
\[ u: {\cal U} \arrow\Teich\] be its pullback.
The map $u$ is a smooth complex analytic
submersion with fiber $(M,I)$ at the point 
$I\in \Teich$ (throughout the paper, $(M,I)$ denotes
a manifold $M$ equipped with a complex structure $I$).
In this paper we use the action of the mapping
class group $\Gamma$ on this
fibration to study the geometry of families of
rational curves on $(M,I)$.

Fix a cohomology class $z\in H^2(M,\Z)$.
Let $\Gamma_z\subset \Gamma$ be the stabilizer of $z$ in
$\Gamma$, and $\Teich_z$ the Teichm\"uller space of all
complex structures $I\in \Teich$ such that
$z$ is of Hodge  type (1,1) on $(M,I)$.

Recall that the second cohomology group
 of a hyperk\"ahler
manifold with maximal holonomy is equipped
with a canonical bilinear symmetric
pairing $q$, called Bogomolov-Beauville-Fujiki (BBF) form. This form is integral but in general not unimodular,
so that it embeds $H_2(M, \Z)$ into $H^2(M,\Q)$ as an overlattice of $H^2(M,\Z)$. It is often
convenient to consider the homology classes of curves as second cohomology classes with rational
coefficients, and we do this throughout the paper. 
Assume that $q(z,z)<0$ and $z$ is represented by a rational curve on some $(M,I)$. It turns out 
that there is a subspace $\Teich_z^{\min}$, which is the same as $\Teich_z$ up to inseparability 
issues, such that for all $I\in \Teich_z^{\min}$, rational curves with cohomology class proportional to $z$ exist on $X=(M,I)$ and 
can be contracted birationally (\ref{_Teich_min_Definition_} and Section \ref{_MBM_loci_Section_}). For $I$ such that the complex manifold $X$ is projective, this is a consequence
of Kawamata base point free theorem
(\ref{_Kawamata_bpf_Theorem_}). For non-algebraic deformations, it follows from the work by Bakker and
Lehn (\cite{_Bakker_Lehn_}) provided that $b_2(M)>5$ (\ref{_existence_of_contra_Theorem_}). 

We are interested in the behaviour of the contraction loci (that is, the exceptional loci of the contraction maps) $Z_I\subset (M,I)$ as $I$ varies in $\Teich_z^{\min}$.
These loci are also obtained as the unions of rational (or all) curves of class proportional to $z$.
The crucial fact is that for any integral class $z\in H_2(M,\Z)$
with $q(z, z) <0$,  the action of the group
$\Gamma_z$ on $\Teich_z$ is also ergodic on each
connected component. Moreover we can classify, in the same way as for $\Gamma$ acting on $\Teich$,
the orbits of the $\Gamma_z$-action on the space $\Teich_z^{\min}$ (\ref{orbits-teichz}).

The subvarieties $Z_I$ are exchanged by the action of
$\Gamma_z$ on $\cal U$ and are thus isomorphic along an orbit, which is often dense.
However, when not in the same orbit, 
they can be very different as complex varieties.
Our main purpose in this paper is to show, under some restrictions, that the $Z_I$ form a trivial family 
in the real analytic category (\ref{_exce_sets_real_analy_trivialized_Theorem_}). 

Note that the real analytic manifolds do not have continuous moduli:
indeed their deformations are controlled by the first cohomology
of the tangent bundle, and higher cohomologies of a coherent sheaf in
real analytic category are always zero (see \cite{_Cartan_}, Th\'eor\`eme 3, for submanifolds of $\R^N$, 
and \cite{GMT}, Theorem 2.7 p. 116, \cite{Nar}, p.931, for the reduction to this case). However, singular
real analytic varieties might have continuous moduli.
The easiest way to see this is to look at configuration $C$ of 4 real
lines in $\R P^2$. If these lines intersect in one point, the 
corresponding tangent cone (which is determined intrinsically 
by the real analytic geometry of the pair $(\R P^2, C)$)
is 4 lines in a vector space. The cross-ratio of these 4 lines 
gives a real analytic invariant of this pair.

Those phenomena are dealt with by Thom-Mather theory.
This theory defines stratified diffeomorphism of real
analytic varieties as a homeomorphism inducing a diffeomorphism
on open strata of a stratification of a manifold by singularities.
Thom and Mather proved that in this category real analytic
varieties have no continuous moduli
(see for example \cite{_Mather:stability_}). 
Later, T. Mostowski and A. Parusi\'nski 
proved that this diffeomorphism is a bi-Lipschitz
equivalence
(\cite[Theorem 1.6]{_Parusinski:Lip_stra_subana_}).
We shall see that the deformation of $Z_I$ and related spaces,
such as the corresponding component of the Barlet space and the
incidence variety, are trivial in stratified
diffeomorphism and in the bi-Lipschitz category
(\ref{_main_incidentn_Theorem_}).

The fact that the family of $Z_I$, as $I$ varies, is locally trivial in the real analytic category even though $Z_I$
can be singular, is related to the fact that $Z_I$ are
contraction loci and follows from results of \cite{_Bakker_Lehn_}.
Bakker and Lehn refer to a concept of ``locally trivial deformation'' 
introduced by H. Flenner and S. Kosarew in
\cite{_Flenner_Kosarew_}
(see also Section \ref{_contraction_BL_Subsection_}).
Unlike its name would suggest, a ``locally trivial deformation''
is not a deformation $\pi:\; {\cal X}\arrow B$ which is equivalent
to the product $F\times U \arrow U$ locally on $B$. Instead,
it is a deformation which is locally trivial {\em locally
  in ${\cal X}$}. 

We show that a locally trivial deformation induces
a trivial deformation in the real analytic category (\ref{_defo_loc_triv_realanal_Proposition_}).
Thus the family of contracted IHS manifolds
constructed in \cite{_Bakker_Lehn_} is trivialized
real analytically, along with the family of the contraction loci.
For other related families, such as the Barlet spaces and incidence spaces
associated with minimal rational curves, 
our techniques (that is, combining an ergodicity theorem with a result of Thom-Mather type) give bi-Lipschitz and 
stratified diffeomorphic trivializations.

\subsection{Teichm\"uller spaces, MBM classes and locally trivial deformations}

The Teichm\"uler space $\Teich_z$
is a smooth, non-Hausdorff manifold, equipped with
a local diffeomorphism to the corresponding period
space $\Perspace_z:=\frac{SO^+(3, b_2-4)}{SO^+(1, b_2-4)\times SO(2)}$ (alternatively, this is just the orthogonal of $z$ in the usual period space $\Perspace$, seen as a subset of a quadric in the projective space $\P H^2(M, \C)$) , which becomes one-to-one if we glue together the inseparable points.
Following E. Markman \cite{_Markman:survey_}, the set of preimages of a point  in $\Perspace_z$ 
(that is, the set of complex structures inseparable from a given $I$)  
is identified
with the set of the K\"ahler chambers in the positive cone
of $H^{1,1}(M,I)$, so that each K\"ahler chamber can be seen as the K\"ahler cone of the corresponding complex structure. The classes relevant for us, those of negative square and represented by a rational curve on some 
$(M,I)$, are the so-called MBM classes 
(Subsection \ref{_MBM_intro_Subsection_} and Section \ref{sectMBM}), i.e. such that the orthogonal complement $z^\bot$
contains one of the walls of these K\"ahler chambers
(\ref{_MBM_walls_Theorem_}). Restricting ourselves to the K\"ahler
chambers adjacent to the hyperplane $z^\bot$, we obtain the space $\Teich_{\pm z}^{\min}\subset \Teich_z$. Note that both spaces are non-Hausdorff even at their general points, since there are always at least two chambers adjacent to a 
given wall. Once $z$ is fixed, $z^\bot$ is co-oriented, and we take the set of the
chambers adjacent to $z^\bot$ on the positive side (that is, $z$ must be positive on the K\"ahler cone).
This last space, separated at its general point, is denoted $\Teich_z^{\min}\subset \Teich_z$.  
This is precisely the space of complex structures
$I\in \Teich_z$ such that a positive multiple of $z$ is represented by an extremal rational curve: indeed, by a result of
Huybrechts and Boucksom, the K\"ahler cone is
characterized as the set of $(1,1)$-classes of positive Beauville-Bogomolov square which are
positive on all rational  curves (\cite{_Huybrechts:basic_,_Huybrechts:erratum_}, \cite{_Boucksom-cone_}).

It follows that the boundary of the K\"ahler cone is a union of a ``round part'' (the boundary of the cone of positive-square classes) and locally polyhedral walls (orthogonals to rational curves) which intersect in locally polyhedral faces of higher codimension.\footnote{See e.g.  \cite{HT}, Proposition 13, for
local finiteness issues.} More generally, if $F$ is a face of the K\"ahler cone of $X=(M,I)$ of codimension $k$ in 
$H^{1,1}_{\R}(X)$, then $F$ is contained in (and has a common open part with) the intersection
of several hyperplanes orthogonal to MBM classes $z_1, \dots, z_k$, where $z_i$ are non-negative on 
the K\"ahler cone. We set $\Teich_F=\bigcap_{i=1}^{k}\Teich_{z_i}^{\min}$: this is the part of 
$\Teich$ where all $z_i$ remain, up to a positive multiple, classes of extremal rational curves.

Thanks to Kawamata base-point-free theorem, it is well-known
that if $X=(M,I)$ is projective, the faces of the K\"ahler
cone can be contracted: there is a projective
birational morphism 
$\phi_F: X \to X'$  sending a curve $C$ to a point iff its
class belongs to the subspace $F^\bot$. Conversely, a
projective birational contraction contracts some extremal
face $F$.

Let $X\arrow X_1$ be a birational contraction of $F$ as above.
In \cite{_Bakker_Lehn_}, Proposition 4.5, Bakker and Lehn prove that any, possibly non-projective, small deformation $X_t$ of $X$ such that all $z_i$ 
remain of type $(1,1)$ on $X_t$ contracts onto a ``locally trivial'' deformation of $X_1$ (in the sense of \cite{_Flenner_Kosarew_}), and that all {\it locally trivial} small deformations of $X_1$ appear in this way\footnote{We shall use the term ``birational contraction'' in the non-algebraic setting too,
meaning ``bimeromorphic contraction''.}.

Bakker and Lehn's result has many applications.
The first application, implicit in \cite{_Bakker_Lehn_}, is the existence of 
bimeromorphic contractions for non-algebraic
hyperk\"ahler manifolds, see \ref{_existence_of_contra_Theorem_}. 


Next, we use the locally trivial deformation of the contracted manifold to produce a real analytic trivialization of the
universal family over $\Teich_F$ preserving the
corresponding contraction locus.

As the simplest examples show, the deformation equivalent contraction loci need not be biholomorphic or bimeromorphic.  Our last aim is to show that the fibers of 
their rational quotient fibrations do. To do this, we prove that the diffeomorphisms of contraction loci as above preserve the rational curves.
This is done by establishing similar triviality results for Barlet spaces and incidence varieties, which we only get in the stratified diffeomorphism category.
We use the following
observation. Let $E\stackrel \phi\arrow  B$ be a proper 
holomorphic  (or even real analytic) map, and
assume that $B$ is obtained as
a union of dense subsets, $B= \bigcup_{\alpha \in {\goth
    I}} B_\alpha$, such that for any index $\alpha$,
all fibers of $\phi$ over $b\in B_\alpha$ are isomorphic. Then
all fibers of $\phi$ are homeomorphic, stratified
diffeomorphic and bi-Lipschitz equivalent.

This observation is based on the classical results
by Thom and Mather (we use the version by Verdier \cite{Ver}, particularly well-adapted to our purposes; the bi-Lipschitz case is
due to Parusi\'nski\footnote{\cite[Theorem 1.6]{_Parusinski:Lip_stra_subana_}; 
see also \cite{_Parusinski:Lip_}, 
\cite{_Parusinski:Lip_stra_}.}) 
These affirm that for any
proper real analytic fibration $E\stackrel \phi\arrow B$,
there exists a stratification of $B$ such that the
restriction of $\phi$ to open strata is locally
trivial in the category of topological spaces
(or in bi-Lipschitz category). Since each $B_\alpha$
in the decomposition $B= \bigcup_{\alpha \in {\goth
    I}} B_\alpha$ intersects the open stratum,
this implies that all fibers of $E\stackrel \phi\arrow B$
are homeomorphic and bi-Lipschitz equivalent.

The dense subsets $B_\alpha$ are in our case provided by the
ergodicity of the mapping class group action. Unfortunately for this argument we have to exclude from consideration the complex structures with maximal Picard number, since their mapping class group orbits are closed.

We state our main results precisely in the subsection \ref{main}, after a brief digression on rational
curves in the next subsection.

\subsection{MBM loci on hyperk\"ahler manifolds}
\label{_MBM_intro_Subsection_}

Let $C\subset M$ be a rational curve on a holomorphic symplectic manifold of dimension $2n$. According to a theorem of Ran \cite{_Ran:deformations_},
the irreducible components of the deformation space of $C$ in $M$ have dimension at least $2n-2$.

\hfill

\definition\label{minimal}
A rational curve $C$ in a holomorphic
symplectic manifold $M$ is called {\bf minimal}
if every component of its deformation space in $M$ has dimension $2n-2$ at $C$.

\hfill

\remark\label{uniruled-subv} In \cite{_AV:MBM_}, Section 4, we have defined and studied minimal rational curves in a maximal irreducible uniruled subvariety $Z\subset X$ as curves of minimal degree with respect to a given K\"ahler class. From the proof of \cite{_AV:MBM_}, Theorem 4.4, Corollary 4.6, one sees that this is equivalent to saying that $C$ deforms in a family of dimension exactly $2n-2$ within $Z$. Indeed the fibers of the rational quotient (also called MRC fibration) of $Z$ are $k$-dimensional, where $k$ is the codimension of $Z$ in $X$, and a simple dimension count shows that if $C$ deforms in a family of dimension greater than $2n-2$, then $C$ deforms with two fixed points, splitting into a union of lower-degree rational curves (``bend-and-break lemma''). Therefore all the results of \cite{_AV:MBM_} about minimal rational curves also apply here.

The dimension of $Z$ can take any value between $n$ and $2n-1$. 
Such a subvariety is always coisotropic, and
the rational quotient fibration is equal to the
coisotropic one (i.e. the kernel of the restriction of the
symplectic form is tangent to the fibers)
(\cite{_AV:MBM_}, Theorem 4.4).

\hfill

The key property of a minimal curve is that such a curve deforms together with its cohomology class $[C]$. More precisely,
any small deformation of $M$ on which $[C]$ is still of type $(1,1)$, contains a deformation of $[C]$ 
(\cite{_AV:MBM_}, Corollary 4.8). Taking closures in the
universal family over $\Teich_z$ gives a submanifold of
$\Teich_z$ of maximal dimension (which does not have to
coincide with $\Teich_z$, as it is not Hausdorff) such that every complex structure in this 
submanifold carries a deformation of $C$; this curve,
however, can degenerate to a reducible curve,
and one cannot in general say much about the
cohomology classes of its components (Markman's example on
K3 surfaces is already enlightening, see
\cite{Markman-pex}, Example 5.3).

\hfill

In \cite{_AV:MBM_}, we have defined and studied the MBM classes: these are classes $z\in H^2(M,\Z)$ such that, up to monodromy and 
birational equivalence, $z^{\bot}$ contains a wall of the
K\"ahler cone\footnote{A
  ``wall'' shall always mean a face of maximal dimension,
  that is $h^{1,1}-1$.}. In other words, $z^{\bot}$ contains a wall of some K\"ahler chamber (see \cite{_Markman:survey_} for the definition of the latter, but it amounts to say that those are monodromy transforms of K\"ahler cones of the birational models of $M$). It is clear that 
the Beauville-Bogomolov square $q(z)$ is then negative; on the other hand, one can characterize MBM classes as negative classes 
such that some rational multiple $\lambda z$ is represented by a rational curve on a deformation of $M$ (\cite{_AV:MBM_}, Theorem 5.11, Corollary 5.14; more precisely, on a deformation with Picard group generated by $z$, $\lambda z$ is represented by a rational curve, and on specializations with larger Picard number this rational curve can break up into a reducible one). For our purposes, it is convenient to extend the notion of MBM on the rational cohomology (or integral homology) classes in an obvious way.


Note that it is apriori possible (though we don't have any
examples) that the same rational curve $C$ is contained in
two maximal irreducible uniruled subvarieties $Z_1$ and
$Z_2$ of $M$, in such a way that the deformations of $C$
lying in $Z_1$ form a $2n-2$-parameter family whereas
those lying in $Z_2$ need more parameters. Such a $C$ is,
by our definition, not minimal, but its generic
deformation in $Z_1$ is.

\hfill

\definition
{\bf An MBM curve} is a minimal curve $C$ such that its class $[C]$ is MBM.

\hfill

\definition
Let $C$ be an MBM curve on a hyperk\"ahler manifold $(M,I)$,
and $B$ an irreducible component of its deformation
space (Chow-Barlet space, well-known to be compact when the ambient manifold is compact K\"ahler) in $M$ containing the parameter point for $C$. An
{\bf MBM locus} of $C$ is the union of all curves
parameterized by $B$.

\hfill

As mentioned in the beginning of this subsection, the MBM
loci are coisotropic subvarieties which can have any
dimension between $n$ and $2n-1$, but the family of
minimal rational curves in an MBM locus always has $2n-2$ parameters (see {\cite{_AV:MBM_}, section 4, Theorem 4.4, Corollary 4.6).

\hfill

\definition\label{full-MBM} Let $z$ be an MBM class in $H^2(M,\Q)$. The {\bf full MBM locus} of $z$ is the union of all MBM curves of cohomology class proportional to $z$ and their degenerations (in other words, the union of all MBM loci for MBM curves of cohomology class proportional to $z$). Similarly, if $F$ is a codimension $k$ face of the K\"ahler cone of $M$, orthogonal to a 
$k$-dimensional subspace $N$ in $H^2(M,\Q)$, we define the full MBM locus for $F$ as the union of MBM loci of MBM classes in $N$.

\hfill

\remark \label{ratconn} If the complex structure on $M$ is in
$\Teich_z^{\min}$, the full MBM locus has only finitely
many irreducible components and is simply the union of all
rational curves of cohomology class proportional to
$z$, and similarly for $F$. This is because (by Kawamata base-point-free theorem in the projective case and by Bakker and Lehn's work in general with as assumption on $b_2$\footnote{In the non-projective case without assumptions on $b_2$,this statement can be shown using the density of complex structures corresponding to projective manifolds, but we shall not need it.}, see section 5) there exists a
birational morphism contracting exactly the curves which have cohomology class in $N$ (that is, orthogonal to $F$)\footnote{By a slight abuse of terminology, we say that ``F can be contracted''.} .
 The number of irreducible components
of an exceptional set of a contraction is finite. One knows that these are {\bf uniruled} (\cite{Kawamata}) and by bend-and-break lemma one finds a minimal rational curve in each. In fact the bend-and-break lemma gives minimal curves in the fibers of the rational quotient and this assures that they are contracted to points, see e.g. \cite{_Bakker_Lehn_}, Prop. 4.11, together with \cite{_AV:MBM_}, proof of Theorem 4.4. These results also show that on a holomorphic symplectic variety the fibers of a contraction map coincide with fibers of the rational quotient of the exceptional locus, in particular the fibers of the contraction map are {\bf rationally connected}.


\hfill

\remark \label{proportional-unequal} Answering a question of the referee, let us mention without giving details of the calculation (for which we refer to our paper \cite{AVX}) that on a given manifold there can exist rational curves with proportional cohomology classes which are both minimal 
(and MBM). Indeed, let $\tau:S\to \P^2$ be the double covering of the projective plane ramified along a sextic, and $X=S^{[6]}$ the Hilbert scheme of length-six subschemes of $S$. 
It is a  classical fact that $H^2(X, \Z)=H^2(S, \Z)\oplus \Z E$ where $E$ is
half of the big diagonal class. Also let us denote by $h\in H^2(S, \Z)$ the class of the inverse image under $\tau$ of a line in $\P^2$. 

The direct sum above is orthogonal with respect to the Beauville-Bogomolov intersection form. Using this form, view the 
classes of curves on $X$ as elements of $H^2(X, \Q)$. Any line bundle of degree 6 on $\tau^{-1}C$ where $C$ is a smooth conic in $\P^2$ has at least a pencil of sections, giving a rational curve $R$ on $X$. In other words, one obtains a rational curve by varying 6-tuples of points on the inverse image of a conic. Moreover this is a minimal rational curve of cohomology class $2h-E$  whose deformations cover a divisor, say $D$, the MBM locus of $R$.

Likewise, since the space of sections of a line bundle of degree 4 on a genus two curve is three-dimensional, we obtain a projective plane in $X=S^{[6]}$ from the inverse image of a line
(by taking those 6-tuples of points on $S$ of which four are on the inverse image and two remaining points are fixed). A line $l$ in this plane has cohomology class $h-E/2$. It is likewise minimal and its deformations cover a codimension-two subvariety $Z$, its MBM locus.

In this particular example $Z$ is a part of $D$. Indeed ``4 points on a line'' is a special case of ``6 points on a conic'': just draw a line through the two remaining points.

\subsection{Main results of this paper}\label{main}

We concentrate on the space $\Teich_z^{\min}\subset \Teich_{\pm z}^{\min}\subset \Teich_z$ described in  
the first subsection. Recall that to construct $\Teich_z^{\min}\subset \Teich_z$, we
first take the complex structures where $z^{\bot}$ actually contains a wall of the K\"ahler cone obtaining $\Teich_{\pm z}$, then
take the ``positive half'' (the complex structures such that $z$ is positive on the K\"ahler classes) of it. 
 On the space $\Teich_z^{\min}$, there is an action of the subgroup of
 the monodromy group preserving $z$, and 
it turns out, thanks to the negativity of $z$, that  
almost all orbits of this action are dense. This allows us
to make conclusions such as the uniform behaviour of subvarieties swept out by curves of class $z$
on the manifolds represented by the points of $\Teich_z^{\min}$.

\hfill

\theorem\label{_main_homeo_Theorem_}
Let $M$ be a hyperk\"ahler manifold of maximal holonomy, $b_2(M)>5$, 
and $z\in H_2(M,\Z)\subset H^2(M,\Q)$ a class of negative 
Beauville-Bogomolov square. Assume that $z$ is represented by a minimal rational curve
in some complex structure $I$ on $M$ (this means that $z$ and the curve are MBM, see \cite{_AV:MBM_}, Section 5; deform to the structure where $z$ generates Picard group and use the deformation invariance).
Let $Z=Z_I\subset (M,I)$ be the full MBM locus of $z$.
Then for all $I,I'\in \Teich_z^{\min}$ there exists a real analytic isomorphism
$h:\; (M,I) \arrow (M,I')$ identifying $Z_I$ and $Z_{I'}$. The same holds for the full MBM locus of any face $F$
of dimension $\geq 3$, as the complex structure varies in $\Teich_F$.

\hfill

For the proof, see \ref{_exce_sets_real_analy_trivialized_Theorem_}.
%
%
In Subsection \ref{_Proofs_Subsection_}
we also prove the following variant of
\ref{_main_homeo_Theorem_}.

\hfill

\theorem\label{_main_incidentn_Theorem_}
In the assumptions of \ref{_main_homeo_Theorem_},
let $B_I$ be the Barlet space of all rational
curves of cohomology class proportional to $z$. Then
the map $h:\; Z_I \arrow Z_{I'}$
can be chosen to send any rational curve $C\in B_I$ to some rational curve $h(C)\in B_{I'}$,
inducing a homeomorphism from $B_I$ to $B_{I'}$, for all complex structures $I, I'\in \Teich_z^{\min}$ except possibly those
with maximal Picard number.

\hfill

This in turn yields another version/strengthening of the theorem.
Recall that a uniruled compact K\"ahler manifold has a so-called rational quotient, or
MRC fibration  (\cite{_Campana:MRC_}, \cite{_KMM:MRC_})
whose fiber at a 
general point $x$ consists of all the points which can be reached from $x$ by a chain of
rational curves. In particular, considering such a fibration on a desingularization of a
component of $Z_I$ gives a rational map $Q:\; Z_I\dasharrow Q_I$. Due to the fact that $Z_I$ 
are contractible the map $Q$ is actually regular and coincides with the contraction itself
(cf. Section 4 of \cite{_AV:MBM_} or Proposition 4.11 of \cite{_Bakker_Lehn_}, and also \ref{ratconn}). 

\hfill

\theorem\label{_main_MRC_Theorem_}
In the assumptions of \ref{_main_homeo_Theorem_}, 
\ref{_main_incidentn_Theorem_}, 
consider the contraction maps with exceptional loci
$Z_I$ and $Z_I'$, and let
$Q:\; Z_I\to Q_I$, $Q':\; Z_{I'}\to Q_{I'}$
denote the restriction of the contraction maps to $Z_I$, $Z_{I'}$.
Then $h:\; Z_I \arrow Z_{I'}$ 
induces a bimeromorphism between the fibers of $Q$ and
$Q'$; it is an isomorphism when these fibers are
normal.

\hfill

\proof
We deduce \ref{_main_MRC_Theorem_} from
\ref{_main_incidentn_Theorem_}
as follows. By \ref{ratconn}, the fibers are rationally connected. Any continuous map of rationally connected varieties mapping
rational curves family  to rational curves is automatically
birational. Indeed the tangent spaces to rational curves
span the holomorphic tangent space of the rationally connected variety at a general point, so that such a map sends holomorphic tangent space to the holomorphic tangent space.
However, a homeomorphism between normal
complex analytic spaces which is holomorphic
on a dense open set is holomorphic everywhere.
This result follows from a version of Riemann
removable singularities theorem, see e. g. 
\cite[Theorem 1.10.3]{_Magnusson:cycle_}. 
\endproof

\hfill

Restricting the diffeomorphism $h$ to the irreducible components of the full MBM locus, we obtain the same statements for MBM loci of curves.










See Subsection \ref{_Proofs_Subsection_} for some other variants of the main theorem.

\hfill

\remark One cannot affirm that the same statements hold along the whole of $\Teich_z$, and this is
false already for K3 surfaces. Indeed a $(-2)$-curve on a K3 surface $X$ can become reducible on a suitable 
deformation $X'$. What we do affirm is that in $\Teich_z$ there is another point, nonseparable from the one
corresponding to $X'$, such that on the corresponding K3 surface $X''$ our curve remains irreducible.
In this two-dimensional case, this easily follows from the description of the decomposition into the K\"ahler chambers in
\cite{_Markman:survey_}; \ref{_main_homeo_Theorem_} allows us to go further in the
higher-dimensional case.









\hfill









\section{Hyperk\"ahler manifolds}


\subsection{Hyperk\"ahler manifolds}

To save space, we omit most of the standard preliminaries on 
hyperk\"ahler and holomorphically symplectic
geometry (see \cite{_Besse:Einst_Manifo_} 
and \cite{_Beauville_}). 
%
%
%
%
%
%
%
%
By a {\bf (simple) hyperk\"ahler}, or {\ irreducible holomorphically symplectic (IHS)} manifold 
we mean a simply-connected compact K\"ahler manifold $M$ such that $H^{2,0}(M)$ is generated by a nowhere degenerate form $\Omega$.
When the context requires, we shall also write $(M,I)$ denoting by $M$ the inderlying differentiable manifold and by $I$ a complex structure on $M$.

%
%
%
%
%
%
%
%
%
%
%

On the second cohomology of a hyperk\"ahler manifold there is an integral quadratic form $q$, called 
{\bf Beauville-Bogomolov-Fujiki (BBF) form}. It has signature $(3,b_2-3)$ and is positive
definite on $\langle \Omega, \bar \Omega, \omega\rangle$,
 where $\omega$ is a K\"ahler form. It is of topological origin and can be defined as follows.

\hfill

\theorem (Fujiki, \cite{_Fujiki:HK_})
Let $\eta\in H^2(M)$, and $\dim M=2n$, where $M$ is
hyperk\"ahler. Then $\int_M \eta^{2n}=cq(\eta,\eta)^n$,
for some primitive integer quadratic form $q$ on $H^2(M,\Z)$,
and $c>0$ a rational number.

%

\subsection{Teichm\"uller spaces and the mapping class
  group}
\label{_Teich_subsection_}

\definition 
Let $M$ be a hyperk\"ahler manifold, and 
$\Diff_0(M)$ the connected component of the unity of its diffeomorphism group
({\bf the group of isotopies}). Denote by $\Comp$
the space of complex structures of K\"ahler type on $M$, and let
$\Teich:=\Comp/\Diff_0(M)$. We call 
it {\bf the Teichm\"uller space} of complex structures on $M$.
It is a complex manifold, possibly non-Hausdorff (more generally for Calabi-Yau manifolds, this statement is essentially contained in 
\cite{_Bogomolov:defo_}; see also \cite{_Catanese:moduli_}).


\hfill

\definition 
Let $\Diff(M)$ be the group of 
diffeomorphisms of $M$. We call $\Gamma=\Diff(M)/\Diff_0(M)$ {\bf the
mapping class group}. 



\hfill

If $M$ is IHS, $\Teich$ modulo the subgroup $K\subset \Gamma$ acting trivially on cohomologies is identified with the marked moduli space, which has finitely many
connected components by a result of Huybrechts
(\cite{_Huybrechts:finiteness_}).  We
consider the subgroup $\Gamma_0$
of the mapping class group which preserves the one
containing the parameter point for our chosen complex
structure.

\hfill



\definition  We call the image of $\Gamma_0$ in
$\Aut{H^2(M,\Z)}$ the {\bf monodromy group},
denoted by $\Mon(M)$.

\hfill

\theorem (\cite{_V:Torelli_}) $\Mon(M)$  is a finite index 
subgroup of the orthogonal lattice $O(H^2(M, \Z), q)$.

\hfill

\remark From now on, to avoid heavy notations, 
we denote by $\Teich$ the connected component of the Teichm\"uller space 
containing the parameter point for our given complex structure, and accordingly write $\Gamma$ instead of $\Gamma_0$.

\subsection{The period map}




\definition 
The map $\Per:\; \Teich \arrow {\Bbb P}H^2(M, \C)$ sending $I$ to the line $H^{2,0}(M,I)$ is 
called {\bf the period map}.

\hfill

\remark 
$\Per$ maps $\Teich$ into an open subset of a 
quadric, defined by
\[
\Perspace:= \{l\in {\Bbb P}H^2(M, \C)\ \ | \ \  q(l,l)=0, q(l, \bar l) >0\}.
\]
It is called {\bf the period space} of $M$.

\hfill

\remark 
One has
\[ \Perspace=\frac{SO^+(b_2-3,3)}{SO(2) \times SO^+(b_2-3,1)}=\Gr_{++}(H^2(M,\R)),
\]

the grassmannian of positive planes in $H^2(M,\R)$ (the sign $+$ in $SO^+$ standing for the connected component 
of the unity).
Indeed, the group $SO^+(H^2(M,\R),q)=SO^+(b_2-3,3)$ acts transitively on
$\Perspace$, and $SO(2) \times SO^+(b_2-3,1)$ is the stabilizer of a point. From a complex line $l$ one obtains a real 
oriented plane by taking its real and imaginary part (in that order).

Bogomolov in \cite{_Bogomolov:defo_} proved that the period map
 is a local diffeomorphism, and Huybrechts has shown the surjectivity in \cite{_Huybrechts:basic_}. The second-named author has obtained the following more precise result in \cite{_V:Torelli_}.

\hfill

\theorem
The points of each connected component of $\Teich$ which have the same image in $\Perspace$ are exactly the non-separable points
(so that the period map is the ``Hausdorff reduction'' of a component of $\Teich$, i.e. becomes an isomorphism once the non-separable points are identified).

\hfill

\definition Let $z$ be a class of negative square in $H^2(M,\Z)$. We call $\Teich_z$ the part of $\Teich$
consisting of all complex structures on $M$ where $z$ is of type $(1,1)$.

\hfill

The following proposition is well-known (see e. g. \cite{_Huybrechts:basic_}, 1.14).

\hfill

\proposition\label{teichz} $\Teich_z=\Per^{-1}(z^{\bot})$, 
where $z^{\bot}$ is the set of points corresponding to
lines orthogonal to $z$ in $\Perspace\subset \P H^2(M,\C)$.

\hfill

On $\Teich_z$, we have a natural action of the stabilizer
of $z$ in $\Gamma$, denoted by $\Gamma_z\subset \Gamma$.

\subsection{Ergodicity of the mapping class group action}

%
%
%

\definition
Let $(M,\mu)$ be a space with a measure,
and $G$ a group acting on $M$ preserving the measure.
This action is {\bf ergodic} if all
$G$-invariant measurable subsets $M'\subset M$
satisfy $\mu(M')=0$ or $\mu(M\backslash M')=0$.

\hfill

It is easy to see that most of the orbits of an ergodic action are dense (the union of non-dense ones has measure zero). A theorem of Moore (see \cite{_Moore:ergodi_}) states that a lattice in a non-compact simple Lie group 
$G$ with finite center acts ergodically on
$G/H$, if $H$ is a non-compact subgroup. Taking $G=SO^+(H^2(M,\R))$ and $H$ the stabilizer of a positive two-plane, we deduce that our mapping class group $\Gamma$ acts ergodically on $\Perspace$: indeed the image of 
$\Gamma$ is of finite index in the orthogonal group of $H^2(M,\Z)$, so it is a lattice in $G$.

%
%
%
%
%








%
%

\hfill

In \cite{_Verbitsky:ergodic_} and  \cite{_Verbitsky:erratum_}, Theorem 2.5, a more precise result has been established using
Ratner theory.

\hfill

\theorem\label{orbits-per} Let $L$ be a integral lattice of signature $(a,b)$ with $a\geq 3, b\geq 1$ and 
$a+b\geq 5$, $V=L\otimes \R$, $\Gamma$ a finite index subgroup in $O(L)$. Then there are three types of orbits
of $\Gamma$-action on $\Gr_{++}(V)$:

1) the orbits of rational planes are closed;

2) the orbits of planes containing no non-zero rational vectors are dense;

3) the orbits of planes containing a single rational line $\langle v\rangle$ are ``intermediate'': each irreducible component of the orbit closure consists of all planes containing a given rational vector, and this vector is 
$\gamma v$ for some $\gamma\in \Gamma$.

\hfill

The theorem applies to the action of the mapping class group $\Gamma$ on $\Perspace$ for an IHS manifold $M$ as soon as $b_2(M)\geq 5$, but also to the following situation: 

\hfill

\corollary\label{orbits-perz} Let $M$ is an IHS manifold with $b_2(M)\geq 5+k$,
$z_1,\dots, z_k\in H^2(M,\Z)$ span a negative subspace, and $\Perspace_{z}=\bigcap z_i^{\bot}\subset \Perspace$ be
the locus of period points of complex structures where each $z_i$ is of type $(1,1)$ (the index $z$ here is the multivector 
$(z_1,\dots z_k)$). Let $\Gamma_z$ be the subgroup of $\Gamma$ fixing 
all the $z_i$. Then there are three types of orbits of $\Gamma_z$ on $\Perspace_{z}$: 

1) the orbits of complex structures with maximal Picard number are closed;

2) the orbits of complex structures with no rational vector in the period plane are dense; 

3) the orbit closure of a complex structure whose period plane contains a single rational line $\langle v\rangle$ is a union of irreducible components, each of them consisting of all planes containing a $\Gamma_z$-translate of $v$.

\hfill

The following observation from \cite{_Verbitsky:erratum_} (Proposition 2.7) shall be useful.




\hfill

\proposition \label{_interme_orbit_not_anal_Proposition_}
In the third case of the above orbit classification, each irreducible component of the orbit closure is a fixed point set
of an antiholomorphic involution (with respect to the natural complex structure obtained by identifying
$\Gr_{++}(V)$ with an open subset of a quadric in $\P V_{\C}$, as in the last section). In particular, it is not
contained in any complex submanifold nor contains any positive-dimensional complex submanifold (even locally).

\section{MBM curves and the K\"ahler cone}
\label{sectMBM}


The notion of MBM classes was introduced in \cite{_AV:MBM_}
and studied further in \cite{_AV:Mor_Kaw_}. We recall the setting and
some results and definitions.

First of all, the BBF form on $H^{1,1}(M,\R)$ has signature $(1, b_2-3)$.
This means that the set $\{\eta \in H^{1,1}(M,\R)\ \ |\ \ (\eta,\eta)>0\}$
has two connected components. The component which contains the
K\"ahler cone $\Kah(M)$ is called {\bf the positive cone},
denoted $\Pos(M)$.

The starting point is the following theorem.

\hfill

\theorem
(Huybrechts
\cite{_Huybrechts:basic_,_Huybrechts:erratum_}, Boucksom
\cite{_Boucksom-cone_})
The K\"ahler cone $\Kah(M)$ is the set of all $\eta\in \Pos(M)$
such that $(\eta, C)>0$ for all rational curves $C$.

\hfill

 Observe that it is sufficient to consider the curves of negative square (as 
only these have orthogonals passing through the interior of the positive cone) and which are moreover
extremal, i.e. such that their cohomology class cannot be decomposed as a sum of classes of other
curves. An extremal curve is minimal in the sense of our \ref{minimal}, though apriori the converse needs 
not be true.














The K\"ahler cone is locally polyhedral in the interior of the positive cone (see e.g. \cite{HT}, Prop.13), with
some round pieces in the boundary, and its walls (that is, codimension-one faces) are supported on the
orthogonal complements to the extremal curves.

The notion of an extremal curve is however not adapted to the defor\-mation\--in\-variant context. 
In order to put the theory in this
context we have defined the MBM (monodromy birationally minimal) classes in \cite{_AV:MBM_}\footnote{G. Mongardi has introduced the notion of wall divisors in \cite{Mon}, the two notions turned out to be equivalent.}
Here we recall several equivalent definitions (we refer to section 5 and more specifically to Theorem 5.16 of \cite{_AV:MBM_} for the equivalence).

\hfill

\definition A negative class $z$ in the image of
$H_2(M,\Z)$ in $H^2(M,\Q)$ is called MBM if $\Teich_z$
contains no twistor curves.

\hfill

{\bf Equivalent definitions:} A negative class $z$ is MBM iff a rational multiple of $z$ is represented by a rational curve in some complex 
structure where the Picard group is generated by $z$ over the rationals. 

Also, $z$ is MBM iff in some complex structure $X=(M,I)$ where $z$ is of type
$(1,1)$, the orthogonal complement $\gamma(z)^{\bot}$  contains a wall of the
K\"ahler cone of a birational model $X'=(M,I')$ of $X$ (this is the original definition from 
\cite{_AV:MBM_} which was at the origin of the terminology). 

Moreover in these two equivalent definitions, ``some'' may be replaced by ''all'' without changing the content.

\hfill

\remark A {\bf wall} always means a face of maximal dimension (that is, $h^{1,1}-1$), a change of terminology from 
\cite{_AV:MBM_} where already ``face'' referred only to faces of maximal dimension unless otherwise specified.

\hfill

\theorem \label{_MBM_walls_Theorem_} (\cite{_AV:MBM_}, Theorem 6.2)
The K\"ahler cone is a connected component of the
complement, in $\Pos(M)$, of the union of 
hyperplanes $z^{\bot}$ where $z$ ranges over MBM classes of type $(1,1)$.

\hfill

\definition (cf. \cite{_Markman:survey_}) 
{\bf The K\"ahler chambers} are the other connected components of this 
complement.

\hfill

Moreover we have the following connection between the K\"ahler chambers and the inseparable points of the
Teichm\"uller space (note that for a fixed deformation type the decomposition of $\Pos$ into the K\"ahler chambers is an invariant
of a period point rather than of the complex structure itself, since it is determined by the position
of $H^{1,1}$ in $H^2(M,\R)$):

\hfill

\theorem (\cite{_Markman:survey_}), theorem 5.16) The points of a fiber of $\Perspace$ over a period point 
are in bijective 
correspondence with the K\"ahler chambers of the decomposition of the positive cone of the corresponding 
Hodge structure. Each chamber is the K\"ahler cone of the corresponding complex structure.

\hfill

\definition \label{_Teich_min_Definition_}
The space $\Teich^{\min}_{\pm z}\subset \Teich_z$ is obtained by removing the complex structures 
where $z^{\bot}$ does not support a wall of the K\"ahler cone. 

\hfill

In other words, at a general point of $\Teich_z$, where the Picard group is generated by $z$ over the rationals,
  $\Teich^{\min}_{\pm z}$ coincides with $\Teich_z$, whereas at special points of $\Teich_z$ where we have other MBM 
classes as well, we remove those complex structures where e.g. $z$ becomes a sum of two effective classes,
and rational curves representing $z$ thus cease to be extremal.

\hfill

Notice that the space $\Teich^{\min}_{\pm z}$ is not separated even at its general point, since $z^{\bot}$ divides the positive cone
in at least two chambers. In order to avoid working with such generically non-separated spaces we divide $\Teich^{\min}_{\pm z}$
in two halves: 

\hfill

\definition The space $\Teich^{\min}_{z}$ is the part of $\Teich^{\min}_{\pm z}$ where $z$ has non-negative intersection with K\"ahler classes (that is, $z$ is pseudo-effective).

\hfill

Now at a general point $\Teich^{\min}_z$ coincides with $\Perspace_z$ (but at special points it is still non-separated).


\section{MBM loci and birational contractions}\label{locifam}
\label{_MBM_loci_Section_}

\subsection{Projective case}

\remark Let $z$ be an MBM class in some complex structure $I\in \Teich_z^{min}$. We have defined the full MBM locus $Z$ of $z$
as the union of subvarieties swept out by minimal rational
curves of cohomology class proportional to $z$. By bend-and-break lemma we can find a minimal rational curve through the general point of any component of $Z$, so $Z$ is the union of all 
rational curves $C$ such that $[C]$ is proportional to $z$.

\hfill

%

These loci are interesting since these are centers of elementary birational 
contractions (Mori contractions). In the projective case this is well-known and follows from Kawamata base-point-freeness
theorem.

\hfill

\theorem \label{_Kawamata_bpf_Theorem_}
(Kawamata BPF theorem, \cite{_Kawamata:Pluricanonical_})
Let $L$ be a nef line bundle on a projective manifold $M$ such that $L^{\otimes a}\otimes{\cal O}(- K_M)$ is big
for some $a$.
Then $L$ is semiample.

\hfill

Recall that a holomorphic line bundle $L$ is {\bf nef} if $c_1(L)$ is in the closure of the K\"ahler cone, and {\bf big}
if the dimension of the space of global sections of its tensor powers has maximal possible growth. For the nef
line bundles this last condition is equivalent to 
$c_1(L)^{\dim M}>0$. A {\bf semiample} line bundle is 
a line bundle $L$ such that $L^{\otimes n}$ is base point free
for some $n$; then the linear system of sections of $L^{\otimes n}$ defines a projective morphism with
connected fibers $\phi: M\to M_0$. The bigness of $L$ implies that $\phi$ is birational. Clearly, for a
curve $C$,  $\phi(C)$ is a point if and only if $L\cdot C=0$.

\hfill

\corollary\label{_faces_bir_contractions_Theorem_} 
Let $M$ be a projective hyperk\"ahler manifold. Then faces $F$ of the K\"ahler cone of $M$, except for the rays 
contained in the boundary of $\Pos(M)$, are in bijective correspondence with birational contractions $\pi:\; M \arrow M_1$, and the exceptional set of $\pi$ is exactly the
full MBM locus of $F$.

\hfill 

\proof First of all, note that the result $M_1$ of a birational contraction is itself projective by the singular version of Huybrechts' criterion  (\cite{BL}, Theorem 6.9).  The face $F$ of the K\"ahler cone is a subset with non-empty interior of the orthogonal complement to some rational cohomogy 
classes (those of extremal rational curves $[C_1], ..., [C_k]$, by Huybrechts-Boucksom description \cite{_Boucksom-cone_}), hence it contains an integral point in its interior when $M$ is projective. 
This point is the Chern class of a nef and big line bundle $L$. The bundle $L$ is semiample since $K_M$ is zero, and hence defines a 
contraction $\pi:\; M \arrow M_1$. Conversely, let $\pi:\; M \arrow M_1$ be a birational contraction and let $L_1$ be an ample bundle on $M_1$. Then $L:=\pi^* L_1$ is a big and
nef line bundle with $c_1(L)\in \bigcap_i [C_i]^\bot$, where $C_i$ are the extremal rational curves contracted by $\pi$ (note that
the contraction loci are uniruled, as one deduces for instance from \cite{Kawamata}, Theorem 1). Hence $\bigcap_i [C_i]^\bot$ is 
a non-empty face.
\endproof

\hfill

\subsection{Non-projective case: locally trivial deformations}
\label{_contraction_BL_Subsection_}



The notion of locally trivial deformations was developed
in  \cite{_Flenner_Kosarew_} and applied to hyperk\"ahler geometry
in \cite{_Bakker_Lehn_}.

\hfill

\definition\label{_locally_triv_defo_Definition_}
Let $\pi:\; {\cal X}\arrow B$ 
be a family of complex varieties. Assume that
any point $x\in {\cal X}$ has a neighbourhood $W$ which
is biholomorphic to a product $F\times U$
such that $\pi\restrict{F\times U}$ is a projection to $U$ ($F$ stands for a neighbourhood of $x$ in the fiber and $U$ for the neighbourhood of its image on the base).
Then $\pi$ is called {\bf a locally trivial deformation},
or {\bf locally trivial deformation in the sense of Flenner-Kosarew}.

\hfill

Let $M$ be a hyperk\"ahler manifold and $f:M\to M_1$ a birational contraction which contracts precisely the curves whose classes 
are in the subspace $N\subset H^2(M, \Q)$. Let $\Def(M)$, $\Def(M_1)$ are local deformation spaces and $\cal X$, ${\cal X}_1$ the
universal families. According to Namikawa \cite{Namikawa}, there is a natural commutative diagram extending $f$:

$$\begin{CD}
{\cal X}  @>\Phi>>  {{\cal X}_1}\\
@VVV              @VVV\\
 {\Def(M)}   @>G>>  {\Def(M_1)}
\end{CD}$$

The fiber of ${\cal X}_1$ over a general point of $\Def(M_1)$ is smooth and the restriction of $\Phi$ to a general fiber $X_t$ of $\cal X$
is an isomorphism, so this diagram in itself does not carry information on contractions. An important advance has been recently made
by Bakker and Lehn.

\hfill

\theorem\label{_loc_triv_BL_Theorem_} (\cite{_Bakker_Lehn_}, Proposition 4.5)
Let $f: M \arrow M_1$ be a birational contraction of 
a projective hyperk\"ahler manifold, with $b_2(M_1)\geq 5$.
Let $\Def^{lt}(M_1)\subset \Def(M_1)$ be the 
subspace parametrizing locally trivial deformations of $M_1$ and 
$\Def(M, f)\subset \Def(M)$ be the subspace of deformations of $M$ 
on which the classes contracted by $f$ remain of type $(1,1)$.
Then the contraction induces an isomorphism
between $\Def(M, f)$ and $\Def^{lt}(M_1)$, so that the small deformations of $M$ preserving the Hodge type of the classes contracted by $f$ contract onto locally trivial small 
deformations of $M_1$.
\endproof

\hfill

Let $M \arrow M_1$ be a birational contraction of
a projective hyperk\"ahler manifold obtained from a face $F$ of its
K\"ahler cone $\Kah(M)$
(\ref{_faces_bir_contractions_Theorem_}). 
Assume that $F$ is supported on (i.e. is a subset with non-empty interior of) the intersection of 
orthogonal complements to linearly independent
MBM classes $h_1, ..., h_n\in H^{1.1}(M)$. By analogy with $\Teich_z^{\min}$, we define the space $\Teich_F$.
Namely $\Teich_F$ is the part of the Teichm\"uller
space of $M$ such that for all $I\in \Teich_F$
the orthogonal complement
$\langle h_1, ... h_n\rangle^\bot$ intersects the
closure of the K\"ahler cone of $(M,I)$ in a face
$F_I$ of the same codimension $n$, and all $h_i$ are positive on the K\"ahler classes. In other words $\Teich_F$ is the intersection
of $\Teich_{h_i}^{\min}$. \ref{teichz} easily implies that $\Teich_F$ is a generically Hausdorff 
manifold equipped with the period map
$\Per:\;  \Teich_F\arrow \Gr_{+,+}(W_F)$, where
$W=\langle h_1, ... h_n\rangle^\bot$. By Torelli
theorem $\Per:\;  \Teich_F\arrow
\Gr_{+,+}(W_F)$ is locally a diffeomorphism.

The following theorem is essentially due to Bakker and Lehn, though not stated in \cite{_Bakker_Lehn_} explicitely.\footnote{Note added in proof: it is in the last version, with a reference to the present paper: Corollary 5.9.}

\hfill

\theorem\label{_existence_of_contra_Theorem_}
Let $I\in \Teich_F$ be a 
complex structure on a hyperk\"ahler manifold $M$. Assume that $\dim_\C \Teich_F\neq 1,2$;
this is equivalent to $\dim_\R F\neq 1,2$.
Then there exists 
a birational map $(M,I) \arrow M_1$ which contracts all
curves with cohomology classes orthogonal to $F$ and only those curves.
Moreover, such a map is uniquely determined by
the space $\langle h_1, ..., h_n\rangle$ used
to defined the face $F$. 

\hfill

\remark If $(M,I)$ is algebraic this is true without the extra assumption on $\dim_\R F$
(\ref{_faces_bir_contractions_Theorem_}).
It would be rather suprising if it were necessary in the non-projective case
but we don't know how to avoid it. The assumption in the version of Bakker and Lehn's paper available to us is 
$b_2(M_1)\geq 4$, so that $\dim \Teich_F=2$ is allowed. This seems to be a misprint as their method of proof needs Verbitsky's description of monodromy orbits, available from dimension three on.

\hfill



{\bf Proof of \ref{_existence_of_contra_Theorem_}.} 
For any algebraic $(M,I)$, $I\in \Teich_F$, the face $F$ is contractible by 
a morphism $f:M\to M_1$. By \ref{_loc_triv_BL_Theorem_}, $F$ remains contractible on small deformations of $M$, say over a small
open neighbourhood $U_I\subset \Teich_F$ of $I$. Let now $(M,I')$ be non-algebraic. 
At this point Bakker and Lehn
 use the ergodicity of the mapping class group action 
as follows. Let $\Gamma_F$ be the subgroup of 
the mapping class group preserving the $h_i$. Standard arguments
imply that $\Gamma_F$ is a lattice in the
Lie group $O(W)$, where $W= {\langle h_1, ..., h_n\rangle}^{\bot}\subset H^2(M, \R)$
(see \cite{_AV:MBM_}).  Then one applies the description of
orbit closures in \ref{orbits-per} and \ref{orbits-perz} to obtain that 
the mapping class group orbit of any non-algebraic complex
structure contains an algebraic one in its closure. Hence such an orbit has a representative
in $U_I\subset \Teich_F$ for $I$ algebraic. For such a representative,
all relevant MBM curves can be contracted by
\ref{_loc_triv_BL_Theorem_}. 
However, all complex structures in the same orbit are isomorphic and 
the isomorphism preserves the classes $h_i$.  \endproof

\hfill

\remark The key ingredient of the orbit closures description is the application of Ratner theory to $\Gamma_F$
action on the $O(W)$-homogeneous  space $\Gr_{++}(W)$.
In order for Ratner theory to be applicable, the connected component $H$ of the stabilizer of a point
$v\in \Gr_{++}(W)$ needs to be generated by unipotents.
The Lie algebra of $H$ is 
isomorphic to $\goth{so}(1, \dim F-1)$, and $H$ is 
generated by unipotents if and only if $\dim_\R F>2$, whence the restrictions on $b_2(M)$.
\endproof

\section{Locally trivial deformations and real analytic geometry}\label{mainproof}

Any deformation of a smooth complex manifold is trivial 
(locally on the base) in real analytic category. This is
most easy to see by constructing an {\bf Ehresmann
  connection} and integrating it to obtain a
flow of diffeomorphisms between the fibers. Recall that an Ehresmann connection 
on a smooth family (i.e. such that $\pi$ is a  submersion with compact fibers) $\pi:\; {\cal X}\arrow B$ is
a splitting of the exact sequence
\begin{equation}\label{_family_tangent_Equation_}
0 \arrow T_{{\cal X}/B} {\cal X} \arrow T{\cal X} \arrow
\pi^* TB \arrow 0,
\end{equation}
where $T_{{\cal X}/B} {\cal X}$ denotes the sheaf of
vector fields tangent to the fibers. It is not hard to
see (by integrating local vector fields lifted from the base via the splitting) that the deformation is trivialized over $B$ if and
only if it admits an Ehresmann connection.
Obstructions to the splitting of
\eqref{_family_tangent_Equation_}
lie in $\Ext^1(\pi^* TB, T_{{\cal X}/B} {\cal X})=
H^1({\cal X}, \Homcal(\pi^* TB, T_{{\cal X}/B} {\cal X}))$.
However, on a real analytic variety higher cohomology
of all coherent sheaves vanishes (\cite{_Cartan_}), hence 
this sequence splits, and one can trivialize the
deformation.

For a singular family
$\pi:\; {\cal X}\arrow B$, the splitting does not always
exist, even in the real analytic category (Subsection \ref{_Teich_hk_Subsection_}).
However, ``locally trivial'' (in the sense of 
Flenner and Kosarew) deformations are trivialized.

\hfill

Throughout this section, the base $B$ is assumed to be smooth.

\hfill

\proposition\label{_defo_loc_triv_realanal_Proposition_}
Let $\pi:\; {\cal X}\arrow B$  be a deformation of compact complex
varieties, which is locally trivial in the sense of
\ref{_locally_triv_defo_Definition_}.
Then the real analytic map $\pi_\R:\; {\cal X}_\R\arrow B_\R$  
underlying $\pi$ defines a family which is trivial
over any sufficiently small open set $U\subset B$.

\hfill

\proof 
By Artin's analytification theorem
(\cite{_Artin:equations_}, Cor. 1.6), it would suffice to trivialize
the family $\pi_\R$ in a formal neighbourhood $\hat F$ of 
$F:= \pi^{-1}(b)$, for all $b\in B$.
Denote by $\hat \pi:\; \hat F \arrow \hat B$ the
corresponding map in the mixed formal-analytic category
(the variety $\hat F$ is analytic along $F$ and
formal in the transversal direction). \footnote{The mixed
  formal-analytic setting is natural for the
  deformation theory of complex analytic varieties,
such as in \cite{_KV:Periods_} or in
\cite{_Barannikov_Kontsevich_}.
The objects of the relevant category are 
complex varieties formally completed
in some directions. To be more rigorous,
an object ${\cal X}$ of this catefory is a pro-scheme
obtained as an inverse limit of 
complex analytic spaces with the same reduction $X$.
The formal deformation space of $X$ 
is obtained as such an inverse limit,
hence it belongs to this category.
If, instead of complex analytic,
we start in the category of algebraic
(Noetherian) schemes, the same approach
gives the usual formal schemes.
The analytification of a formal
deformation is a complex analytic
space containing $X$ as a closed
complex analytic subvariety, with the
formal completion along $X$ identified
with ${\cal X}$.}

Locally in ${\cal X}$, the complex family $\pi$ is a
product. The local-in-${\cal X}$ trivialization of $\pi$ 
defines a \v Cech cocycle $w\in H^1(F, \Aut_F(\hat F))$
where $\Aut_F(\hat F))$ is the group sheaf of automorphisms
of $\hat F$ trivial on $F\subset \hat F$ and commuting with the projection to
$B$.
The sheaf $\Aut_F(\hat F))$ can be obtained as a limit
of sheaves of automorphisms of infinitesimal neighbourhood
$F_k \subset \hat F$ of order $k$. Therefore,
$w\in H^1(F, \Aut(\hat F))$ vanishes whenever 
its finite order representatives $w_k\in H^1(F, \Aut_F(F_k))$ vanish. The Lie groups
$\Aut_F(F_k)$ are nilpotent, and 
fit into exact sequences
\[
0 \arrow V_k\arrow \Aut_F(F_k) \arrow \Aut_F(F_{k-1}) \arrow 0
\]
where $V_k$ is a sheaf of abelian unipotent
groups, that is, a coherent sheaf. In the corresponding
exact sequence of first cohomology
\[ H^1(V_k) \arrow H^1(\Aut_F(F_k)) \arrow H^1(\Aut_F(F_{k-1}))
\]
all terms vanish, which can be shown by induction. Indeed,
$\Aut_F(F_1)$ is trivial because the automorphisms commute
with the projection to $B$. On the other hand, higher cohomology of any coherent
sheaf on a real analytic variety vanishes (\cite{_Cartan_}, Th\'eor\`eme 3, completed by \cite{GMT} Theorem 2.7, \cite{Nar} p. 931). 
We obtain that the group sheaf $\Aut_F(F_k)$ is filtered by normal subgroups
with coherent subquotients, hence has vanishing cohomology.
\endproof

\hfill

In the sequel, a
``vector field'' on a singular variety $S$ is understood 
as a section of the sheaf $(\Omega^1_S)^*$ (dual to the
K\"ahler differentials). By the universal property
of K\"ahler differentials, $(\Omega^1_S)^*$
is the sheaf of derivations from $\calo_S$ to itself.

\hfill

\proposition\label{_trivialize_pair_Proposition_}
Let $\pi:\; {\cal X}\arrow B$  be a deformation of complex
varieties, which is locally trivial in the sense of
Flenner-Kosarew (\ref{_locally_triv_defo_Definition_}), and
 $\sigma:\; \tilde {\cal X}\arrow {\cal X}$ 
a simultaneous resolution of singularities. Assume that any vector field on 
the smooth part of ${\cal X}$
can be lifted to a vector field on $\tilde {\cal X}$.
 Then the family
$f=\sigma\circ\pi:\; \tilde {\cal X}\arrow B$
admits a real analytic Ehresmann connection such that the
corresponding flow of diffeomorphisms preserves the
exceptional variety ${\cal Z}$ of $\sigma$, and moreover 
does so fibrewise over its image in ${\cal X}$.

\hfill

\pstep We start by showing that it suffices to prove
existence of an Ehresmann connection preserving ${\cal Z}$
locally in ${\tilde {\cal X}}$. An Ehresmann connection in 
$f:\;  \tilde {\cal X}\arrow B$ is the same as a
splitting of the exact sequence
%
\begin{equation}\label{_family_tangent_for_tilde_X_Equation_}
0 \arrow T_{\tilde {\cal X}/B} \tilde {\cal X} \arrow T\tilde {\cal X} \arrow
f^* TB \arrow 0,
\end{equation}
Therefore, a difference between two Ehresmann
connections is a section of 
$\Homcal(f^* TB, T_{\tilde {\cal X}/B} \tilde {\cal X})$.
Consider the natural pairing
\[
 \Psi:\; J_{\cal Z}/J^2_{\cal Z} \otimes T \tilde {\cal X}\arrow \calo_{\cal Z}
\]
obtained if we identify vector fields with derivations and
take a derivation of $\alpha \in J_{\cal Z}$ evaluating it on ${\cal Z}$.
Clearly, a diffeomorphism associated with 
a vector field $v$ preserves ${\cal Z}$ if and only if
$\Psi(v, \cdot)=0$. This gives a coherent sheaf $\ker \Psi$ denoted by
$T_{\cal Z} \tilde {\cal X}\subset T \tilde {\cal X}$.
This is the sheaf of vector fields preserving ${\cal Z}$.
Now, if we have found an Ehresmann connection preserving ${\cal Z}$
locally in ${\cal X}$, the corresponding \v Cech cocycle $w$ (of differences on intersections) takes
values in 
\[ 
  \Ext^1(f^* TB, T_{\tilde {\cal X}/B} \tilde {\cal X}
  \cap T_{\cal Z} \tilde {\cal X})= H^1((f^* TB)^*\otimes
  (T_{\tilde {\cal X}/B} \tilde {\cal X} \cap T_{\cal Z}
  \tilde {\cal X}));
\]
this group vanishes because cohomology of any coherent sheaf on a 
real analytic variety vanish (\cite{_Cartan_}). Therefore,
the connections constructed locally on ${\tilde{\cal X}}$ give rise to a global one.

\hfill

{\bf Step 2:} 
In Step 1, we reduced \ref{_trivialize_pair_Proposition_}
to a statement which is local on ${\cal X}$. 
Since locally in ${\cal X}$ we have
 ${\cal X}=B \times F$, we can assume that the family
$\pi:\; {\cal X}\arrow B$ is trivial, and ${\cal X}=B \times F$.
This gives a natural embedding $\pi^*TB \hookrightarrow T{\cal X}$.
Replacing $B$ by an open ball if necessary, we fix the
coordinate vector fields
$\underline \zeta_1, ..., \underline\zeta_n\in H^0(TB)$. 
Using the embedding $\pi^*TB \hookrightarrow T{\cal X}$,
we obtain holomorphic vector fields
vector fields $\zeta_1, ..., \zeta_n$ on ${\cal X}$
which can be integrated to diffeomorphisms
${\cal V}^i_t$. These diffeomorphisms 
are coordinate translations along
$B$ in the decomposition  ${\cal X}=B \times F$.

The vector fields $\zeta_i$ can be lifted to holomorphic
vector fields on 
the simultaneous resolution $\tilde {\cal  X}$, 
by assumptions of  \ref{_trivialize_pair_Proposition_}.
Denote the corresponding 
holomorphic diffeomorphism flows on $\tilde {\cal X}$
by $\tilde {\cal V}^i_t$. 
These diffeomorphism flows commute with the
projection $\sigma:\; \tilde {\cal X}\arrow {\cal X}$,
because $\sigma \circ {\cal V}^i_t= \tilde {\cal
  V}^i_t\circ \sigma$ 
at the general point of $\tilde {\cal X}$.
Therefore, the 
diffeomorphism flows $\tilde {\cal V}^i_t$
preserve ${\cal Z}$, and the corresponding
vector fields give a splitting of
\eqref{_family_tangent_for_tilde_X_Equation_}.
\endproof

\hfill

To apply \ref{_trivialize_pair_Proposition_} to holomorphic
symplectic varieties, we use the following lemma.

\hfill

\lemma
Let $\sigma:\;\tilde{\cal X}\arrow {\cal X}$ be a simultaneous birational contraction in a 
family of holomorphic symplectic manifolds over a ball $B$. Then any vector field 
on the smooth part of ${\cal X}$ can be extended to a vector field on $\tilde{\cal X}$.

\hfill

\proof
Notice that the manifold $\tilde{\cal X}$ has trivial canonical bundle, so that 
both on $\tilde{\cal X}$ and on the smooth part of ${\cal X}$, vector fields are identified with
differentials of degree $\dim {\cal X}-1$. The differentials
extend by \cite{KS}, Cor. 1.8.
\endproof

\hfill

Comparing this lemma with \ref{_trivialize_pair_Proposition_},
we obtain the real analytic Ehresmann connection preserving
the exceptional sets of birational contractions:

\hfill

\theorem\label{_exce_sets_real_analy_trivialized_Theorem_}
Let $M$ be a hyperk\"ahler manifold, and $\pi:\; M \arrow M_1$
a birational contraction associated with a face $F$ of the
K\"ahler cone of $M$. Assume that $b_2(M_1)\geq 5$, and
consider the universal family ${\cal U} \arrow \Teich_F$ of 
over the Teichm\"uller space $\Teich_F$, and the 
corresponding universal family of birational contractions 
${\cal U} \stackrel\sigma \arrow {\cal U}_1 \arrow \Teich_F$  constructed
by Bakker and Lehn (see \ref{_contraction_BL_Subsection_}).
Then the family ${\cal U} \arrow \Teich_F$ admits a 
real analytic trivialization which preserves the fiberwise
exceptional sets of the contraction $\sigma$.
\endproof

\hfill

Our main \ref{_main_homeo_Theorem_} obviously follows.


\section{Applications of Thom-Mather-Verdier theory to the
  families of MBM loci}


We shall now prove a weaker form of \ref{_main_homeo_Theorem_} for e.g. the family of Barlet spaces. In our previous paper \cite{AV:loci}
which the current one supersedes, this method has been applied to the initial family of MBM loci $Z_I\subset M_I$ for which a stronger result has just been
obtained using Bakker-Lehn's theorem. Our old method is based on two ingredients which apply in a great generality, thus permitting to obtain 
a weaker result for essentially any family related to the geometry of rational curves on an IHSM.



\hfill

One ingredient is the work by Verdier on the Whitney stratification and Thom-Mather theory in the 
complex analytic context \cite{Ver}. It implies that the members of any proper complex analytic family (i.e. the fibers $X_y$
of a proper morphism of countable at infinity complex analytic spaces $f:X\rightarrow Y$) are homeomorphic and 
stratified diffeomorphic (with respect to a strong Witney stratification, \cite{Ver}, 2.1) over a complement 
 to a union of closed analytic subvarieties.


The other ingredient is the description of the orbits of the monodromy action on $\Teich_z^{\min}$, or more generally $\Teich_F$, which is the same as the one 
for the period space but the proofs are somewhat more technical (\ref{orbits-teichz}). This description allows to send, by an element 
of the mapping class group, a point on such 
a subvariety (along which a topological/stratified differentiable degeneration in a family over $\Teich_z^{\min}$ or 
$\Teich_F$ is supposed to happen) into a small neighbourhood of a general point. As the mapping class group acts by diffeomorphisms, this proves that the degeneration actually does not happen,
unless the Picard number at that point is maximal (in this case the mapping class group orbit is closed so the argument does not work).

We now give the details of the argument sketched above, restricting for simplicity of notation to the families over $\Teich_z^{\min}$ (but the argument is the same over 
$\Teich_F$ which is the intersection of several $\Teich_z^{\min}$).

\subsection{Mapping class group action on $\Teich_z^{\min}$}

The group $\Gamma_z\subset \Gamma$ obviously acts on $\Teich_z^{\min}$. Indeed the action of any 
$\gamma\in \Gamma$ is just the transport of the complex structure; if $z^{\bot}$ contains a wall of the
K\"ahler cone in a complex structure $I$, then so does $\gamma z$ in the complex structure $\gamma I$.
Notice that the same remark applies to rational curves: $\gamma C$ is a rational curve in the structure
$\gamma I$ and the minimality is preserved. So the full MBM locus $Z\subset X=(M,I)$ of $z$
is sent by an element of $\Gamma_z$  to the full MBM locus $Z_{\gamma I} \subset X'=(M, \gamma I)$. 

It turns out that the results on the mapping class group action on $\Perspace$ ``lift'' to those on the action on 
$\Teich$, but if we want to work on a subspace where $z$ remains of type $(1,1)$ this has to be $\Teich_z^{\min}$
rather than $\Teich_z$.

The following theorem from \cite{_Verbitsky:ergodic_}, \cite{_Verbitsky:erratum_} strengthens \ref{orbits-per}.

\hfill

\theorem\label{orbits-teich} 
Assume $b_2(M)\geq 5$. Let $\Gamma$ denote the mapping class group.
Then there are three types of $\Gamma$-orbits on $\Teich$: 
closed (where the period planes are rational,
thus the complex structures have maximal Picard number), dense (where the period planes contain no rational
vectors), and such that each irreducible component of the closure is formed by points whose period planes contain a fixed rational vector $\gamma v$ (where $\gamma\in \Gamma$ and $v$ generates the unique rational line in the period plane).
In the last case, no neighbourhood of a point $c$ in the orbit closure $C_v$ is contained in a proper complex subvariety of $\Teich$.

\hfill

The argument in \cite{_Verbitsky:erratum_} proves that the closure of $\Gamma I$ is the inverse image of the closure of $\Gamma \Per (I)$ when the Picard group of $I$ is not maximal. The key idea is to replace 
$\Per: \Teich \to \Perspace$ by a $\Gamma$-equivariant embedding $\Teich_K\to \Per_K$, where $\Teich_K$ consists of pairs $(I,\omega)$ where $I\in \Teich$ and $\omega\in \Kah(I)$ is of square one, and $\Per_K$ consists of pairs $(\Per(I), \omega)$, 
$\Per(I)\in \Perspace$, $\omega\in \Pos(I)$ of square one (note that the positive cone is an 
invariant of the 
period point). Since the points of $\Teich$ correspond to pairs 
$(\Per(I), k)$ where $\Per(I)\in \Perspace$ and $k$ is a K\"ahler chamber of $\Pos(I)$, it suffices to prove that the closure of the orbit of $(I, \Kah(I))$ as a set contains the orbit of $(\Per(I), \Pos(I))$.
The space $\Per_K$ being homogeneous, one uses Ratner theory to do this.

\hfill

The analogue of \ref{orbits-teich} in our setting is as follows. The interesting feature is that it is
$\Teich_z^{\min}$ rather than $\Teich_z$ which replaces $\Teich$.

\hfill

\theorem\label{orbits-teichz}
Assume $b_2(M)>5$. Let $z\in H^2(M,\Z)$ be an MBM class and $\Gamma_z$
the subgroup of the mapping class group consisting of all elements whose action on the second cohomology fixes $z$. Then $\Gamma_z$
acts on $\Teich_z^{\min}$ ergodically,
and there are the same three types of orbits
of this action as in \ref{orbits-teich}.

\hfill

\proof It proceeds along the same lines as in \cite{_Verbitsky:erratum_}. We introduce the
spaces $\Per_{K,z}$ consisting of pairs 
\[
\{(\Per(I),
\omega)| I\in \Teich_z, \omega \in \Pos(I)\cap z^{\bot}, q(\omega, \omega)=1\}
\]  
and $\Teich_{K,z}$ consisting of pairs $(I, \omega)$ where $I\in
\Teich_z^{\min}$, 
and $\omega$ of square 1 
belongs to the wall of $\Kah(I)$ given by $z^{\bot}$.
We denote such a wall by  $\Kah(I)_z$, though of course 
its elements are not K\"ahler forms on $I$, but rather
nef limits of those. Since the complex
structures in $\Teich_z^{\min}$ which 
have the same period point are in one-to-one
correspondence with the walls of the K\"ahler chambers in
which the other MBM classes partition $z^{\bot}$, $\Teich_{K,z}$
again embeds naturally in $\Per_{K,z}$. We fix a complex
structure $I$ with non-maximal Picard number. We need to
prove that the closure of the $\Gamma_z$-orbit of 
the subset $(I,\Kah(I)_z)$ contains the orbit of
$(\Per(I), \Pos(I)\cap z^{\bot})$. This is done exactly in
the same way as in \cite{_Verbitsky:erratum_}, proof of theorem 3.1. The key idea in \cite{_Verbitsky:erratum_}
is as follows: the non-maximality of Picard number means that the subspace generated by the integral classes in 
$H^{1,1}_I(M,\R)$ has non-zero orthogonal complement $T$. Consider a three-dimensional subspace $W\subset H^{1,1}_I(M,\R)$ 
of signature $(1,2)$ such that $W\cap T\neq 0$. Then the intersection of the K\"ahler cone with 
$W$ has a ``round part'' (\cite{_Verbitsky:erratum_}, Proposition 3.4) and therefore contains horocycles (\cite{_Verbitsky:erratum_}, subsection 3.3). In our context, the space $T$ is clearly contained in $z^{\bot}$, meaning that for $W\subset z^{\bot}$
of signature $(1,2)$ intersecting $T$, $\Kah(I)_z\cap W$ contains horocycles. As in \cite{_Verbitsky:erratum_}, Proposition 3.5 and the following paragraph, we deduce from  Ratner's
orbit closure theorem that the closure of the projection
of such a horocycle to $\Per_{K,z}/\Gamma_z$ is large, containing an $SO(H^{1,1}(I)\cap z^{\bot})$-orbit, which is the
projection of  $\Pos(I)\cap z^{\bot}$.
\endproof

\subsection{Stratification}

Consider the family of hyperk\"ahler mamifolds ${\cal X}$ over
$\Teich_z^{\min}$ ($X_I=(M,I)$ over each point $I$, one can introduce it as a pullback of the universal family from (\cite{_Markman:universal_}).
Throughout this paper we have been interested in the family ${\cal Z}\subset {\cal
  X}$ with the fiber over $I\in \Teich_z^{\min}$ obtained
as the full MBM locus of $z$ on the complex manifold
$X=(M,I)$. This family can be constructed
by taking the image of the evaluation map for the union of the components of the relative Barlet space corresponding to
cohomology classes proportional to $z$ and dominating 
$\Teich_z^{min}$. 

Another family we shall consider is the dominating part of the relative Barlet space itself: in such a way we
obtain a family ${\cal B}$ over $\Teich_z^{\min}$, and we call $B_I$ the fiber over $I$ (the Barlet space of minimal 
rational curves in classes proportional to $z$, it has compact components because $X_I$ are compact K\"ahler). 

Finally, there is the incidence variety ${\cal J}\subset {\cal X}\times_{\Teich_z^{\min}} {\cal B}$. 
As $\Teich_z^{\min}$ is not Hausdorff,
we shall, whenever necessary, restrict all families to a 
small neighbourhood $U$ of some point $x$, or to a small compact 
$K$ within $U$, and denote by ${\cal X}_U$, ${\cal Z}_U$ etc. 
the restrictions of these families. 

Whitney \cite{W}  introduced stratifications of analytic varieties by singularity type. Recall that a stratification of $Y$ is a finite
filtration by closed subsets $Z_i$ of dimension $i$ such that each difference $Z_i-Z_{i-1}$ is a manifold, in general non-connected; the strata are its connected components. To use stratifications in practice, one needs some ```glueing conditions'' of technical nature, the so-called Whitney's A and B conditions, or Verdier W condition
(these are equivalent in the complex analytic case). A Whitney stratification is a locally finite stratification satisfying those conditions. Verdier (\cite{Ver}, Th\'eor\`eme 2.2) proved that
a complex analytic space, countable at infinity, admits a Whitney stratification by complex analytic strata, and moreover such that a given closed analytic subset is a union of strata.

Thom-Mather theory uses stratifications to prove the stratified differentiable local triviality of a family over an open subset of the base, for example via the following  {\bf first isotopy
lemma} (see \cite{_Mather:stability_} for a detailed but accessible account).

\hfill

\lemma (\cite{_Mather:stability_}, 
Proposition 11.1)  Let $f: Y\to B$ be a smooth mapping of smooth manifolds and $W$ a closed subset of $Y$ admitting 
Whitney stratification, such that $f:W\to B$ is proper. If the restriction of $f$ to each stratum of $W$ 
is a submersion, then $W$ is locally trivial over $B$ (topologically and stratified differentiably).

\hfill

In the situation of this lemma, one says that $f$ is a ``controlled submersion'', or ``transverse to the stratification''. Verdier (\cite{Ver}, Th\'eor\`eme 3.3) proves a very general result to the effect that a proper morphism of complex analytic spaces $f:X\to Y$ is transverse to a stratification of the source over the complement to a closed analytic subset (more generally, if $f$ is proper in restriction to a closed subset $A$ with a stratification, then $f$ is transverse to the stratification over a dense open subset of $f(A)$). 

The first isotopy lemma holds for a morphism of complex analytic spaces $f:X\to Y$, with non-singular $Y$, proper over a stratified closed subset $W$. Th\'eor\`eme 4.14 of \cite{Ver} is then an analogue of the first isotopy lemma and implies local topological and stratified differentiable triviality of a proper morphism of 
complex spaces over a complement to a closed analytic subset (\cite{Ver}, Corollaire 5.1, formulated as an example of what one obtains in the algebraic case, but the proof remains valid with the properness hypothesis in the
analytic case).

The results of $\cite{Ver}$ applied to our situation yield the following

\hfill

\lemma\label{loctriv} The family ${\cal Z}_U$ is locally topologically and stratified differentiably trivial 
over a complement to a (lower-dimensional) analytic subset, and so is ${\cal B}_U$. Moreover 
$({\cal X}\times_{\Teich_z^{\min}}{\cal B})_U$ admits a trivialization preserving the incidence subset ${\cal J}$.

\hfill

{\it Proof:} This follows from the above recollection of Verdier's results, keeping in mind that the family 
${\cal B}_U$ is proper over $U$ by the compactness of the cycle spaces for compact K\"ahler manifolds, and 
so, by the same reason, is ${\cal Z}_U$. Finally, the preservation of ${\cal J}$ is a consequence of it being a 
union of strata for a suitable Whitney stratification (\cite{Ver}, Th\'eor\`eme 2.2).

\hfill

\remark Concerning ${\cal Z}_U$, this result is of course weaker than the one already proved using 
contractibility. Unfortunately, this other method does not seem to apply to Barlet and incidence spaces.

\subsection{Proof of \ref{_main_incidentn_Theorem_} and closing remarks} 
\label{_Proofs_Subsection_} 

\ref{_main_incidentn_Theorem_} is a consequence of the following fact.

\hfill

\theorem If $b_2(M)>5$, the families ${\cal B}$ and ${\cal J}\subset {\cal X}\times_{\Teich_z^{\min}}{\cal B}$ are topologically (and stratified-differentiably) trivial over the whole  
${\Teich_z^{\min}}$, with a possible exception of points corresponding to the complex structures with maximal Picard number.

\hfill

\proof
We know by \ref{loctriv} that this is the case 
over the complement to a union (possibly countable,
but finite in a neighbourhood of any point
in the base) of proper analytic subsets $P=\cup_iP_i \subset \Teich_z^{\min}$.
First we pick a point $x\in \Teich_z^{\min}$ which is not
in $P$ and whose $\Gamma_z$-orbit is dense.
Then $x$ has a neighbourhood $U_x$ over which all fibers ${\cal B}_b$, $b\in U_x$ are homeomorphic. Moreover the union 
$\bigcup_{\gamma\in \Gamma_z}\gamma(U_x)$ is a dense open subset of $\Teich_z^{\min}$ and all fibers ${\cal B}_b$
over this union are homeomorphic.

Take another point $x'\in \Teich_z^{\min}$ (which now can be in $P$) with dense $\Gamma_z$-orbit (i.e. ``ergodic''). Then the orbit of $x'$ hits 
$\bigcup_{\gamma\in \Gamma}\gamma(U_x)$. But $\Gamma_z$ is a subgroup of the mapping class group and its action is just the transport of the complex structure. Therefore rational curves in a complex structure $I$ and in 
$\gamma(I)$ correspond via $\gamma$, and so do the MBM loci, Barlet spaces, incidence varieties. So ${\cal B}_{x'}$ is homeomorphic to ${\cal B}_b$ for 
$b\in U_x$ and no degeneration happens at $x'$.

Now take $y\in \Teich_z^{\min}$ such that the corresponding complex structure is not ergodic but does not 
have maximal Picard number either
(``the intermediate orbit'' of \ref{orbits-per}
and \ref{orbits-perz}). 
If ${\cal B}_{y}$ is not homeomorphic to ${\cal B}_b$ for 
$b\not\in P$, the orbit of $y$ should remain in $P$ and so must the orbit closure. Each irreducible component
of the orbit closure must be
contained in an irreducible component of $P$, but this is an analytic subvariety.
However, the closure of an intermediate orbit is not
contained in a proper analytic subset, even locally 
(\ref{_interme_orbit_not_anal_Proposition_}), so this is impossible.

The proof for the stratified diffeomorphic case is exactly the same, and 
also the same arguments apply to ${\cal J}$.
\endproof

\hfill

In  \ref{_main_incidentn_Theorem_}
we prove that the fibers of natural
families associated with rational curves 
are homeomorphic and stratified diffeomorphic.
However, there is a version of the Thom-Mather theory
which gives bi-Lipschitz equivalence of the fibers
over open strata of Thom-Mather stratification
(\cite[Theorem 1.6]{_Parusinski:Lip_stra_subana_};
see also \cite{_Parusinski:Lip_}, \cite{_Parusinski:Lip_stra_}).
Then the same arguments as above
prove that the homeomorphisms constructed
in \ref{_main_incidentn_Theorem_} are bi-Lipschitz.

\hfill

{\bf Acknowledgements:}
We are grateful to Fedor Bogomolov for pointing out a potential error
in an earlier version of this work, and to Jean-Pierre
Demailly, Patrick Popescu, Lev Birbrair and Daniel Barlet for useful
discussions. We are especially grateful to
Fabrizio Catanese who explained to us the basics
of Thom-Mather theory and gave the relevant 
reference, and to A. Rapagnetta and the anonymous referee
of the superseded version of the paper for bringing Bakker
and Lehn's paper to our attention and insisting on its
importance for our subject. Much gratitude
is due to Grigori Papayanov for insightful comments
and the reference in Mathoverflow \cite{_Mathoverflow:Papayanov_}.
The referee of the present version has done a considerable work pointing out our many inaccuracies, we thank him/her very much. Remark \ref{proportional-unequal} is inspired by a conversation with Emanuele Macri.

\hfill

{\scriptsize

{\small
\noindent {\sc Ekaterina Amerik\\
{\sc Laboratory of Algebraic Geometry,\\
National Research University HSE,\\
Department of Mathematics, 6 Usacheva Str. Moscow, Russia,}\\
\tt  Ekaterina.Amerik@gmail.com}, also: \\
{\sc Universit\'e Paris-11,\\
Laboratoire de Math\'ematiques,\\
Campus d'Orsay, B\^atiment 425, 91405 Orsay, France}

\hfill

\noindent {\sc Misha Verbitsky\\
{\sc Instituto Nacional de Matem\'atica Pura e
              Aplicada (IMPA) \\ Estrada Dona Castorina, 110\\
Jardim Bot\^anico, CEP 22460-320\\
Rio de Janeiro, RJ - Brasil }\\
also:\\
{\sc Laboratory of Algebraic Geometry,\\
National Research University HSE,\\
Department of Mathematics, 6 Usacheva Str. Moscow, Russia,}\\
\tt  verbit@impa.br}.
 }
}

\end{document}